# Global optimization test problems based on random field composition

R. Sala*, N. Baldanzini, M. Pierini

Dipartimento di Ingegneria Industriale, Università degli Studi di Firenze, Via di Santa Marta, 3 - 50139 Firenze Italy
* rsala.unifi@gmail.com

**Abstract**

The development and identification of effective optimization algorithms for non-convex real-world problems is a challenge in global optimization. Because theoretical performance analysis is difficult, and problems based on models of real-world systems are often computationally expensive, several artificial performance test problems and test function generators have been proposed for empirical comparative assessment and analysis of metaheuristic optimization algorithms. These test problems however often lack the complex function structures and forthcoming difficulties that can appear in real-world problems. This communication presents a method to systematically build test problems with various types and degrees of difficulty. By weighted composition of parameterized random fields, challenging test functions with tunable function features such as, variance contribution distribution, interaction order, and nonlinearity can be constructed. The method is described, and its applicability to optimization performance analysis is described by means of a few basic examples. The method aims to set a step forward in the systematic generation of global optimization test problems, which could lead to a better understanding of the performance of optimization algorithms on problem types with particular characteristics. On request an introductive MATLAB implementation of a test function generator based on the presented method is available.

*Keywords:* global optimization, performance test problems, test function, artificial landscapes, problem features, metaheuristics, random fields, variable interactions

## 1      Introduction

Many tasks in science and engineering can be formulated as optimization problems. A wide variety of optimization algorithms that attempt to solve such problems has been introduced and reviewed in the literature [1, 2]. There are several different criteria that can be used to classify optimization problems [1, 3, 4], and sometimes even the classification of a well-known problem can be a quite difficult problem itself [5]. Common classes of problem types are for example: NP-hard (Nondeterministic Polynomial time) for the level of "difficulty", and non-convexity for the "shape" of the objective function landscape. Many industrial and real life problems belong to the category of non-convex problems, and their solution can be difficult or at least "expensive" [6] in terms of effort.

To deal with the difficulties in the large class of non-convex problems, many different metaheuristic algorithms (MHAs) have been developed since the second half of the last century (see [7, 8] for reviews of such algorithms). Even though some of these algorithms might not always result in the global optimum, often no better alternatives to deal with these problem types are available. It remains however a challenge, to select the most efficient algorithm (i.e. which uses fewer function evaluations) among all the available ones to optimize a specific problem type. Since theoretical performance analysis is difficult and still only restricted to simple MHAs and particular problem types, "empirical" performance analysis based on numerical experiments are common practice in the field of MHAs [9]. Many different test functions and problems have been proposed, and several function collections were published [10, 11, 12, 3, 13]. Furthermore, practical testing environments containing test problem collections and benchmarking frameworks have been programmed [14, 15, 16, 17], and are widely used in the optimization community. However in recent works [18-20] many of the common test functions have been criticized because they are not enough challenging and do not represent the features and types of difficulty[1] as real-world problems.

In [18] the topic of test function generators for assessing the performance of MHAs on multimodal functions was discussed. Authors highlighted, that many of the currently available test functions in the specialized literature are too simple, and show regularities such as: symmetry, uniform spacing of optima, and centered optima. Since such regularities can be easily exploited by algorithm designers (see also [20]), those functions are unrepresentative testing environments. Several strategies to generate more complicated and realistic test functions were proposed [21-24]. However, to the knowledge of the authors, none of available approaches deal in a systematic way with function features such as: higher order design variable interactions, variance contribution distributions, and nonlinearity. An exception is a parallel work by the authors [26], on representative surrogate problems that involved the generation of test functions with particular variance contribution distributions in an application-oriented context. In order to analyze the influence of particular function features on the optimization algorithm performance, it would be of interest to parameterize test problems w.r.t. feature characteristics, and to investigate multiple distinct problem instances of similar problem types.

Recently the issues regarding variable interactions and non-separability have been addressed in a survey [27] of MHAs in Large-Scale Global Optimization problems (LSGO). The survey concluded that more effort is needed to develop decomposition-based optimization methods, such as Cooperative Coevolution methods [28], with high performance on non-separable and separable subcomponents and that their performance on imbalanced problems should be further investigated. Although the specific application to LSGO problems is not further addressed in the present paper, the proposed method can generate test problems with properties of high relevance to this topic (unbalanced variance contribution distributions, and degree of non-separability). Moreover, the method is scalable in order to create test functions for large-scale problems.

The proposed method enables the systematic construction of test problems with varying structures with parameterized variance contribution distributions, higher order interactions, and heterogeneous modality. Furthermore multiple problem instances with the same problem specifications can be generated, which facilitates the statistical assessment of MHA performance on different instances of a problem type. Besides as stand-alone test functions, the fields or functions can also be added to existing test functions to enrich their complexity and increase the level of difficulty.

---

[1] When "difficulty" or "hardness" is averaged over all possible search or optimization algorithms no problems are intrinsically harder than others [25]. However, for a particular optimization algorithm, some problem classes can be more difficult than others.

The aim of this communication is to introduce a concept: to devise structured functions that are based on the superposition of random fields and to demonstrate that they can provide useful insights as test functions in the field of global optimization. The concept of the method is described in section 2. In section 3 the presented method is applied to investigate the influence of changes in function features on the optimization performance of a genetic algorithm. The last sections comprise discussion, outlook, and conclusions.

## 2    Description of the method

In this section a concept for test function generation based on the composition of random fields is presented. This section is divided in three parts:

1. The description of A basic multidimensional discrete random field generator capable to produce parameterized fields with higher order interactions
2. The description of a "smoothing" method to obtain continuous and differentiable fields, by means of weighting functions.
3. A description of several composition techniques to create structured fields by combining different types of basic random fields

*A basic Multidimensional Discrete Random Field (MDRF) generator*

Random fields are of interest in various branches of Mathematics, Physics and Engineering. A random field is a stochastic process taking values in a Euclidean space [29]. Elementary discrete random fields (DRF) can be interpreted as a list of "random" numbers with the indices mapped onto an *n*-dimensional space. The general idea presented in this communication is to compose a number of discrete random fields of different spatial resolutions and dimensionality, in order to construct fields with particular structures. Such fields can serve as test functions of highly variable difficulty in terms of spatial nonlinearity, variance contribution distributions, and higher order interactions.

To model computationally affordable fields which can possess higher order interactions we describe a MDRF generator. This generator function (referred to as operator **H**) takes a multidimensional vector $x$ of floating point values from the unit hypercube domain as an input, and maps it to a value from a given set $S$ with a discrete probability distribution and a computational type (e.g. binary integer or float) of choice:

$$y = \mathbf{H}(x_1, x_2, x_d, \ldots x_n) \text{ where } x_d \in [0,1] \text{ for d=1,2,…,}n \text{ and } y \in S \qquad (1)$$

To explain the concept of the implementation, a related discretized version of this idea can be defined as:

$$y = \mathbf{A}(j_1, j_2, j_d, \ldots j_n) \text{ where } j_d \in \mathbb{N} \text{ and } j_d \leq r_d \text{ and } y \in T \qquad (2)$$

Operator **A** can be interpreted as high-dimensional random "array" *A* with indices *j* of which each index $j_d$ is bounded by the maximum array size $r_d$ for dimension *d*. In expression 2, $T$ is a finite set of successive integers pointing to the distinct elements of the set $S$.

The concept chosen for the MDRF generator algorithm is to compute and reproduce the pseudo-random values of the high-dimensional arrays "on the fly" instead of storing a potentially huge passive map in the computer memory[2]. Another alternative interpretation would be to consider operator **A** as a pseudo-random number generator with a high-dimensional vector as its generating seed.

The ideas in equations 1 and 2 can be combined to establish a parameterized MDRF operator: $y = \mathbf{H}(x, r, \varphi, S)$ where the operator on input variables $x$ is parameterized with respect to the discretization resolution $r$, discretization offset $\varphi$ and codomain set $S$.

The concept of the algorithm can be explained by the following steps:

1. Addition of an optional offset or shift $\varphi_d$ to design variable vector $x_d$: $\hat{x}_d = x_d + \varphi_d$
2. Discretize the resulting floating point input variables $\hat{x}_d$ to integers with respect to the resolution $r_d$ or corresponding array size : $j_d = \text{ROUND}(\hat{x}_d, r_d)$, where $j_d$ represents an index vector.
3. Map the resulting index vector $\boldsymbol{j}$ to an integer index $i$ of **T** by a Pseudo-Random Mapping (PRM) : $i = \text{PRM}(\boldsymbol{j})$
4. Return the element $y \in S$ to which the resulting integer index $i$ pointed: $y = S_i$

The pseudo-randomness and higher order interactions of the resulting discrete field are introduced by the PRM (step 3). The PRM can be achieved by using a Pseudo-Random Number Generator (PRNG), with a multivariate random seed mechanism. Depending on the program implementation such an approach would easily enable the use of discrete random fields with a total array size $m$ of $10^{1000}$ or more ($m = r^n$ where n is the dimensionality of the problem). In the following sections the parameterized implementation of the MDRF generator is denoted by $\mathbf{H}(x, r, \varphi, S)$ or primitives thereof (when parameters not of interest in a particular context are omitted for better readability).

*Obtaining continuous and differentiable random fields*

The algorithm presented in the previous section can generate multidimensional discrete random fields, with specific probability density distributions. Fig. *1* shows on the left a two-dimensional example of such a discrete field, with an array size or spatial resolution of five intervals per dimension in the domain.

---

[2] Such a passive map would be very memory intensive since the required memory scales with the number of elements $m = \prod_{d=1}^{n} r_d$ or $m = r^n$ for a uniform resolution $r$ and field dimension $n$. which already becomes problematic at modest resolutions and problem dimensions. A discrete field array of resolution $r = 10$ and dimension $n = 12$ would already require 8 TB (terabyte) of memory when each element takes 8 bit of storage.

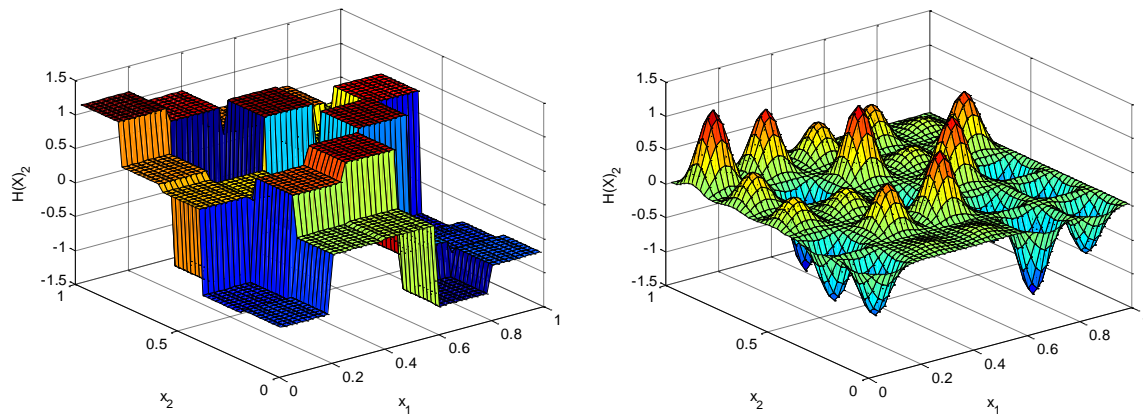

**Fig. 1** *Left* Examples of a discrete random field in 2 dimensions. *Right* and the corresponding continuous field after smoothing

Using a multiplicative weighting function as expressed in equation 3 these discrete random fields can be transformed to "smooth" differentiable continuous random fields.

$$\mathbf{C}(\mathbf{x}, \mathbf{r}, \boldsymbol{\varphi}, p) = \left( \prod_{d=1}^{n} \left( \left( 1 - e^{\cos(2\pi(r_d x_d + \varphi_d)) - 1} \right) * \left( \frac{1}{1 - e^{-2}} \right) \right) \right)^{p} \tag{3}$$

Where $\boldsymbol{\varphi}$ denotes the vector of the (phase) shift of the discrete field and the $\mathbf{C}(\ )$ operator such that also fields with non-zero values along the domain boundaries can be constructed ($\varphi_d \in [0,1]$). With parameter $p$ the shape of the function can be adjusted. This weighting function has the properties, that at each location where the discrete random field is not differentiable in one or more directions, the value and corresponding derivative of the weighting function are equal zero for all values of $p > 0$ [3]. Fig. **1** (right) shows a "smoothened" version of the discrete field using this method. The multiplicative composition of operators $\mathbf{H}()$ and $\mathbf{C}()$ is abbreviated as $\widetilde{\mathbf{H}}()$. Although the application of the $\mathbf{C}(\ )$ operator or weighting function, is technically not smoothing, we will refer to it as smoothing to avoid misunderstandings with the weighting factors introduced later. This "smoothing" operator works to generate continuous fields from the discrete fields in arbitrary dimensions, however it should be noted that in this context as an effect of high dimensionality the integral over the product of the smoothing or weighting function can become vanishingly small with a rate that depends on the choice of exponent $p$. This effect is similar to the decreasing relative volume of the *n*-dimensional hypersphere with respect to the volume of the unit hypercube for high dimensions. The smoothing operator also affects the probability density distribution of the resulting field w.r.t. the original discrete field distribution. For high-dimensional spaces these effects can be controlled by choosing exponent $p$ sufficiently small, which decreases the "smoothness" but still enforces continuity and differentiability. Alternatively other smoothing approaches could be considered.

---

[3] Although the above statement is true for all $p > 0$ the author recommends to use as a rule of thumb $p \geq 1/d$, since for very small values of $p$ the smoothness vanishes.

*Random Field Composition (RFC) based test functions*

The application of the "bare bones" discrete random fields generated by the algorithm in the previously described sections, as optimization test functions is of little practical interest because of the primitive problem structure. The message of this section is however that compositions of such fields of different and heterogeneous resolutions, dimensions and codomain distributions can provide test functions with interesting problem structures.

Continuous fields with different spatial resolution can be created, and compositions can be made by for example multiplication or by weighted addition such as for example:

$$\widetilde{\mathbf{H}}^{comp}(x) = \sum_{k=1}^{m} w_k * \widetilde{\mathbf{H}}_k(x, \boldsymbol{r}_k, \boldsymbol{\varphi}_k) \tag{4}$$

where $\boldsymbol{r}_k$ and $\boldsymbol{\varphi}_k$ (both in bold) denote the vectors with the array size and shifts for each dimension of composition fields *k*. A graphical example of such a weighted field summation in 2d is displayed in **Fig. 2**. For clearness of visualization, only two fields of low resolution were added, but additions of many fields, with higher and distinct resolutions are possible.

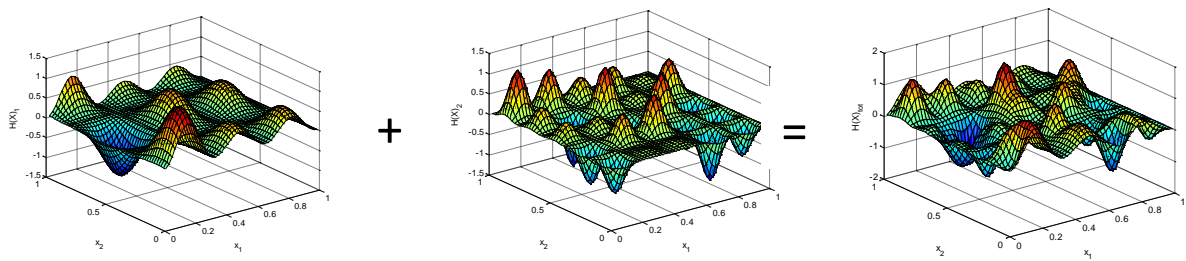

**Fig. 2** Graphic example of the addition of 2 smoothened discrete fields of different resolutions and the resulting composition

According to the Sobol-Hoeffding decomposition [30, 31], it is possible to decompose a vector valued function $\mathbf{F}(x)$, into unique summands of increasing dimensions when the summands are orthogonal.

$$\mathbf{F}(x) = f_0 + \sum_{i=1}^{n} f_i(x_i) + \sum_{1 \leq i < j \leq n} f_{i,j}(x_i, x_j) + \cdots + \sum_{1 \leq i < j < \cdots \leq n} f_{i,j,\ldots,n}(x_1, x_2, \ldots, x_n) \tag{5}$$

In expression 5 multi-index notation is used[4] [32]. The different terms of the summation, refer the unique interaction terms of all possible combinations of variable subsets. In the field of sensitivity analysis, this idea is commonly used to decompose the variance contribution of the function output variance w.r.t the individual summands, to identify the variance contribution of variables and interactions of variable subsets. Such results are commonly expressed by means of quantitative importance measures named Sobol indices *D*, or sensitivity indices, that can be defined as in eq. 6 [32]:

$$D_{i,j,\ldots,n} = \mathrm{Var}\left(f_{i,j,\ldots,n}(x_1, x_2, \ldots, x_n)\right)/\mathrm{Var}(\mathbf{F}(x)) \text{ with } 1 \leq i < j < \ldots \leq n \tag{6}$$

---

[4]Explanation to the multi-index notation: The expression $\sum_{1 \leq i < j \leq n} f_{i,j}(x_i, x_j)$ indicates a sum over all function decomposition terms with two variables for which $1 \leq i < j \leq n$. This applies similarly to all pairs of higher order interactions $f_{i,j,\ldots,n}$.

This expresses the variance contribution of a subset of variables as: the ratio of the variance of the corresponding term[5] from the Sobol-Hoeffding decomposition to the variance of the total function. This concept can be used to quantify (additive) function (output) separability w.r.t. its (input) design variables.

In this context, we apply these ideas in order to construct functions with predefined variance contribution distributions of the first and higher order interaction terms by specifying the weights in equations 4 and 7 accordingly. Besides composition of fields over a fixed set of dimensions or design variables (such as the example in figure 2), also fields over different variable subsets can be composed in order to generate fields targeting variance contributions by specific interaction terms. Equation 7 shows how a random field can be composed of weighted sums of random fields over variable subsets.

$$\widetilde{\mathbf{H}}_{tot}^{comp}(x) = w_0 + \sum_{i=1}^{n} w_i * \widetilde{\mathbf{H}}_i^{comp}(x_i) + \sum_{1 \leq i < j \leq n}^{i.n} w_{i,j} * \widetilde{\mathbf{H}}_{i,j}^{comp}(x_i, x_j) + \cdots + w_{i,j,\ldots,n} * \widetilde{\mathbf{H}}_{i,j,\ldots,n}^{comp}(x_1, x_2, \ldots, x_n) \quad (7)$$

Particular variance contribution distributions over a selection of subsets can be achieved by applying the weights according to one's needs. Each of the subset fields $\widetilde{\mathbf{H}}_{i,j,\ldots,n}^{comp}$ can themselves also be composed of a summation of fields over the corresponding variable subset (see equation 4). A notable point is however, the uniqueness and orthogonally of the summands. In general, the random vectors or fields generated for variables or variable subsets are not necessarily orthogonal to each other. For subfields of high resolution, or high dimensionality the lower-order interaction effects will average out and will become approximately orthogonal w.r.t lower-order effects of other fields in the composition. For lower resolutions and dimensionality such separable effects cannot be neglected and have to be accounted for by for example an *a posteriori* sensitivity analysis on the final composed function, or a covariance/correlation coefficient analysis between the composition summands.

A test function $\mathbf{Z}(x)$ based on the presented concept of parameterized MDRF composition can then be described by the general expression:

$$\mathbf{Z}(x) = \widetilde{\mathbf{H}}_{\overline{w},\overline{r},\overline{\varphi},\overline{S}}^{comp}(x) = \widetilde{H}(x, \overline{w}, \overline{r}, \overline{\varphi}, \overline{S}) \quad (8)$$

where now the composition parameters $\overline{w}, \overline{r}, \overline{\varphi}, \overline{S}$ indicate arrays/structures containing the all the parameter vectors of the composed fields.

The 2d "landscapes" from the previous visualization examples, are not really any more spectacular than landscapes of existing test functions. The novelty of the method lies in the parameterization of the function structure (with respect to variance contribution distributions, function modality and higher order interactions) combined with the straightforward scalability to create high-dimensional test problems.

## 3  examples: RFC based optimization algorithm performance analysis

For a few example problems, the isolated effects of some function features on the optimization performance of a genetic algorithm are demonstrated. The optimization algorithm used is a simple genetic algorithm (SGA) from the publicly available Genetic Algorithm Toolbox for MATLAB developed by Chipperfield et al. [33]. For the investigations a population size of 1000, combined with default settings were used. The optimization performance is measured in the number of function evaluations N that is required to find a solution within $\varepsilon$ of the best known solution ($\varepsilon = 10^{-3}$).

---

[5]The variance for the terms expression 6, w.r.t. the corresponding sub domain in the unit hypercube can be expressed as: $\text{Var}\left(f_{i,j,\ldots,n}(x_1, x_2, \ldots, x_n)\right) = \int f^2_{i,j,\ldots,n}(x_1, x_2, \ldots, x_n) \, dx_i \ldots dx_n$.

*Variance contribution distributions*

The first example shows the influence of the different distributions of the first order sensitivity indices or variance contribution. For a given instance of the first order term related fields, the weights $w_k$ can be optimized such that a particular distribution for the sensitivity indices for the first order terms $D_i$ can be obtained. The level of difficulty in terms of SGA function evaluations to convergence is compared for different distributions in figure 3. The target first order sensitivity index distributions of the small 10 dimensional example problem are chosen according to:

$$\sigma_i = \frac{(i/n)^k}{\sum_{i=1}^{n}\left(\frac{i}{n}\right)^k} \text{ and } k \geq 0 \tag{9}$$

Such that different types of distributions (uniform for $k = 0$, linear for $k = 1$, and skewed for $k > 1$) can be obtained. This example only adresses first order sensitivity index distributions. However, also the distribution of higher order effects are expected to influence the optimization. The example in figure 4 shows that for increasing values of $k$, and decreasing effective dimension, the problem gets significantly easier to solve for the selected optimization algorithm.

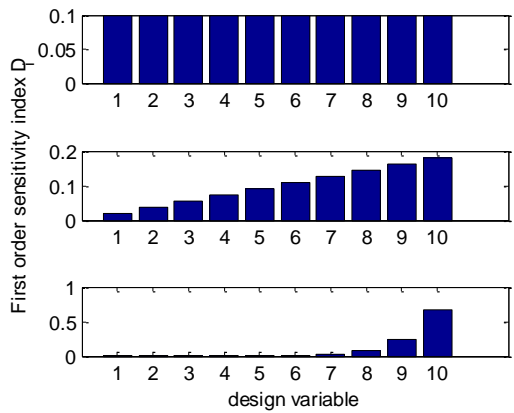
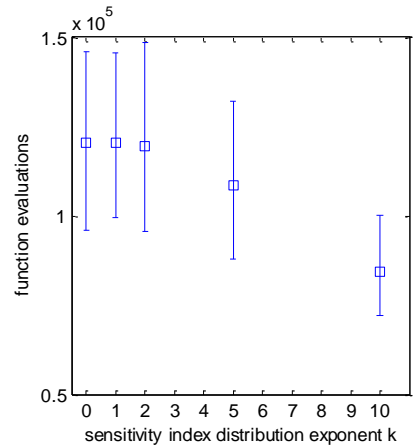

**Fig. 3** Variance contribution distribution examples for *k* values: 0 (top); 1 (center); 10 (bottom)

**Fig. 4** For an increasing first order sensitivity index distribution exponent *k*, the average number of required function evaluations decreases. (DRF settings: d=10, r=20)

*Variable interaction order*

The second example shows the effect of increasing design variable interactions on the level of difficulty of the problem expressed in the number of function evaluations required to converge. In equation 7 the different term types represent different interaction orders. Each interaction order $q$ adds $\binom{n}{q}$ interaction terms i.e. discrete random fields with $r^q$ degrees of freedom. For each interaction term the corresponding weights ($w_{i,j,k,l,m}$) are chosen such that the variance contribution of the corresponding weighted field is: $\frac{1}{Q\binom{n}{q}}$ where $Q$ is the maximum order of interaction of the problem. For the small example problem of dimension 5 with interaction orders up to 5, figure 5 shows as expected that for increasing interaction order both the mean level of difficulty as well as the variance increase.

*Multimodality and DRF resolution*

The third example, shows how the optimization algorithm efficiency scales w.r.t. the chosen base resolution of the discrete random field. This is done for a small 5 dimensional smoothened random field containing only first order interactions. In this example the resolution is homogeneous, and the DRF values are taken from a uniform distribution. Thus, the "resolution" directly affects the multimodality of the resulting test function. The modality is further dependent on the choice of the targeted distribution **S** for the discrete field. As expected the required number of function evaluations increases with increasing resolution.

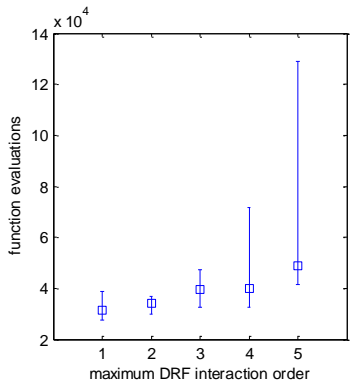 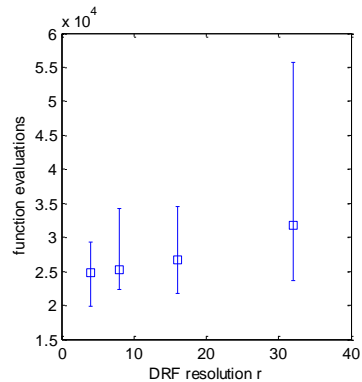 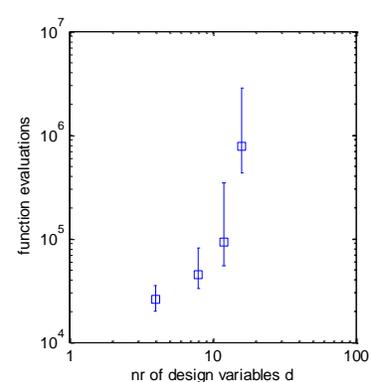

**Fig. 5** The average number of required function evaluations N, and the 20% and 80% percentiles, for an increasing interaction order (DRF settings: d=5 and r=5)

**Fig. 6** The required number of function evaluations for increasing DRF resolution (DRF settings: d=5)

**Fig. 7** The required number of function evaluations for increasing DRF dimension (DRF settings: r=10, interaction order ~d)

*Dimensional scalability*

A common algorithm performance discrimination criterion is the scaling behavior of optimization algorithms on test function instances with increasing dimensionality. In the presented approach parameterization of the dimensionality of the problem is straightforward. A fourth example shows the number or function evaluations required by the genetic algorithm until the convergence criterion is met. Figure 7 shows super-polynomial scaling. In this example, a single smoothened DRF was created such that the interaction order increases proportionally to the number of design variables. The example shows that quite easily difficult problems can be created that require many function evaluations to solve. It is however also possible to limit the interaction order of the problems so that the scaling will be less strong.

*Deception*

Also, test functions containing the established function features such "deception" can be constructed. The general idea behind deceptive functions is a global trend or spatially "larger" function basins that "distract" the optimization algorithm away from the smaller basin of the true global optimum. This effect can be achieved in a statistical sense, by the composition of smoothened lower resolution fields, combined with high resolution fields of which the probability density function is such that: the "amplitude" of these fields is negligible, except at a single or few places where the magnitude of the amplitude dominates the amplitudes of the lower resolution smoothened fields.

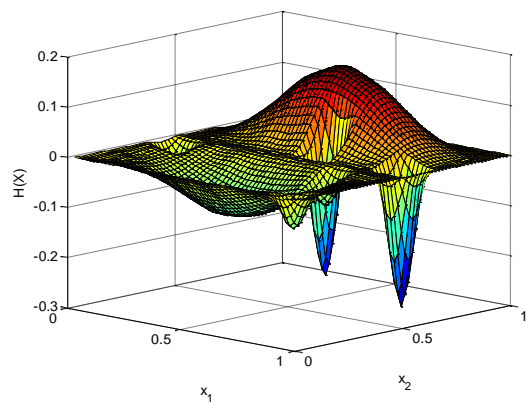

**Fig. 8** Example of a DRF based deceptive function

## 4    Discussion and outlook

These few examples in the previous section, briefly demonstrate the potential of the concept of RFC based test problems to systematically construct optimization problems. Several different function properties were parameterized and the influence on the optimization performance was investigated. In most examples isolated properties are assessed, but the RFC based function generation concept covers a vast function space containing many combinations of different "landscape" properties.

When analyzing the performance of optimization algorithms on different problems, one is interested in which algorithms or algorithm operators can exploit certain properties of the problem structure. The novelty and additional value of the proposed method over other methods is the tunable level of difficulty regarding modality, variance contribution, and design variable interaction order.

A drawback of the presented method is that for compositions of these discrete random field based test functions the global optimum is not necessarily known (depending on the composition type). Probability based estimations can however be made. Alternatively brute force evaluation could be used to explore the details of the generated search space. Although a priori knowledge of the global optimum is desired for test functions, such knowledge is often also not available for real-world problems, and thus this disadvantage could therefore also be seen as a feature of realism.

Another property of the test function generating concept presented is that according to one's wishes it is possible to generate test functions with a huge amount of descriptive parameters. This property is typically classified as unwanted, because of its potential complexity usage. Referring to the concept of Kolmogorov complexity [35], one could however also argue that simple test functions with few parameters are intrinsically very unlikely to represent the difficulties that can arise in highly specified complex computational models, and therefore are strongly limited in their scope. In this context, the authors view the possibility of creating highly parameterized test problems as a necessary feature to specify problems of highly structured complexity.

As also expressed in [18] the development of more challenging test functions could lead to the development of more robust and effective optimization algorithms. In the examples relatively small MHA performance evaluation studies are presented on a single algorithm only. In future work these investigations will be extended to more optimization algorithms, and many more combinations of function features.

The large function space that the method covers, naturally leads to the need for the specification of subproblem classes (parameter ranges) and to define standardized test function instances for algorithm benchmarking. Such a classification could also contribute to compare this method with the many available "anecdotal" test functions, and it could place those functions in a more general function feature context. Other topics for further related research involve the construction of constrained and multi-objective test functions using a similar systematic approach. The flexible parameterization, of the dimensionality, and higher order interactions between groups of design variables is also of interest for the generation of test problems for LSGO. The authors hope that this communication will motivate others to use and extend the concepts presented, in endeavors to analyze and develop useful MHAs.

## 5      Summary and Conclusions

An algorithm is described that can be used to generate and high-dimensional pseudo-random discrete fields of heterogeneous resolution. The resulting discrete fields can be combined with suitable weighting functions to obtain continuous and differentiable fields. Based on composition of these random fields functions with interesting characteristics can be constructed. These test functions can be parametrized w.r.t. a variety of problem features such as: variance contribution distributions, variable interactions, nonlinearity, and dimension.

The resulting method enables systematic performance analysis of optimization algorithms on problems with such function features. The influence of different several function features on the performance of a simple genetic algorithm was investigated by means of a few examples. The method covers a large space of optimization problems, and many other investigations on MHA performance analysis are realizable by means of this method. The primary objective of this communication is to motivate to exploration of the possibilities of RFC based test functions, as a tool for optimization algorithm analysis. Moreover, to stimulate the development and identification of efficient global optimization algorithms for particular problem types. This will hopefully also lead to improved optimization efficiency in real-world optimization problems.

**Acknowledgements**

The presented method is a generalized derivative of strategies developed for efficient optimization on structural simulations in vehicle design in the scope of the GRESIMO project funded by the European Commission under the 7[th] framework program under grant agreement 290050.